\documentclass{article}
\usepackage{amsmath}
\font\tengoth=eufm10
\font\sevengoth=eufm7
\font\fivegoth=eufm5
\newfam\gothfam
\textfont\gothfam=\tengoth \scriptfont\gothfam=\sevengoth
\scriptscriptfont\gothfam=\fivegoth

\def\blacksquare{\hbox to .60em{\vrule width .60em height .60em}}

  \font\bb=msbm10 
\catcode`\@=12  

\def\hb {\hfil \break}
\def\n {\vskip 0.2cm \noindent }
\def\scirc{\,{\raise 0.8pt\hbox{$\scriptstyle\circ$}}\,}
\def\ins{\,{\raise 0.2cm \hbox{ $\scriptstyle \circ$}}\,}

\def  \é{\'e}
\def\è{\`e}
\def\à{\`a}
\def\ù{\`u}
\def\ç{\c c}
\date{}

\begin{document}



 \centerline {\Large\bf Rank   of ordinary webs in codimension one }
 
 \bigskip
 
  \centerline {\Large\bf  An effective method } 

 \bigskip
\rightline {  \bf  J.P. Dufour and D. Lehmann}
 
 \bigskip


 \bigskip
 \n  {\bf Abstract :}
 
 \bigskip
 
 {We  are interested by  holomorphic $d$-webs $W$  of codimension one in a complex $n$-dimensional \hb manifold $M$. If they are \emph{ordinary}, i.e. if they satisfy to some condition of genericity (whose precise definition is recalled below),   we proved in [CL] that their rank $\rho(W)$ is upper-bounded by a certain  number $\pi'(n,d)\ \bigl($which, for $n\geq 3$, is stictly smaller than the 
  Castelnuovo-Chern's bound $\pi(n,d)\bigr)$. 
  
  In fact, denoting by $c(n,h)$ the dimension  of the space of homogeneous polynomials of degree $h$ with $n$ unknowns, and by $h_0$ the integer such that 
  $$c(n,h_0-1)<d\leq c(n,h_0),$$
 $\pi'(n,d)$ is just the first number of a decreasing sequence of positive integers 
 $$\pi'(n,d)=\rho_{h_0-2}\geq \rho_{h_0-1}\geq \cdots\geq \rho_{h}\geq \rho_{h+1}\geq\cdots\geq \rho_{\infty}=\rho(W)\geq 0 $$
 becoming stationary equal to $\rho(W)$ after a finite number of steps.  This sequence  is an interesting  invariant of the web, refining the data of the only rank.
 
 The method is effective : theoretically,  we can compute    $\rho_h$ for any given $h$ ; and, as soon as two consecutive such numbers are equal ($\rho_h=\rho_{h+1}, \ h\geq h_0-2$), we can  construct a holomorphic vector bundle 
 $R_h\to M$ of rank $\rho_h$,  equipped with a tautological holomorphic connection $\nabla^h$ whose curvature $K^h$ vanishes iff the above sequence is stationary from there. Thus, we may stop the process at the first step where the curvature vanishes.

 Examples will be given.

 
 
 
 

  \vskip 1cm
  
   \n  {\bf Contents :}
   
    \bigskip
  
  1- Introduction
  
  \bigskip
  2- Computation of $R_h$
  
    \bigskip
  3- The connections $\nabla^h$
  
    \bigskip
  4- Algorithm
  
    \bigskip
  5- Examples
  
    \vskip 4cm
 \n {\bf Keywords :} ordinary webs, abelian relation, rank, connection, curvature.
 
\n {\bf AMS classification  :} 53A60 (14C21, 53C05, 14H45)

 \pagebreak
 \section{Introduction } 

Recall that a totally decomposable\footnote{More generally,  the web is defined    by one foliation   $\tilde{\cal F}$   on a covering space    $\tilde M\to M$ with $d$ sheets. On an open set   $U$ of  $M$ on which this covering space is trivial,  the data of $\tilde{\cal F}$ is equivalent to  that of the   $d$'s   $ {\cal F}_i$ which are the  projections of  $\tilde{\cal F}$. All  
  fiber bundles and   connections studied below    may be  defined globally.  But all computations being done locally,   we shall recall the definitions only in the  case of a    totally decomposable web.}  holomorphic \emph{$d$-web} of codimension one without  singularity on a complex $n$-dimensional manifold$M$ is   defined by the data of  $d$  holomorphic regular  foliations ${\cal F}_i$ of  codimension one on $M$, $(1\leq i\leq d)$, any one of them being  transverse to each other at any point.

 We assume $d>n$ and the web to be  at least in weak general position\footnote{The web is said to be in weak (resp. strong) general position if, at any point $m$, there exists at least  $n$ of the foliations among the  $d$'s, whose tangent spaces  at $m$ are in general position  (resp. if any  family  of  $n$ foliations  among the  $d$'s have this   property).}.
 
An    \emph{  abelian relation}  on an open set  ouvert $U$ (assumed to be connected and simply connected) of $M$ is then the data of a   family $(F_i)_i$ of holomorphic fèunctions on  $U$, $1\leq i\leq d$, such that 
 
 - for any $i$, $F_i$ is a first  integral of  ${\cal F}_i$ (maybe with   singularities), 
 
 - the sum  $\sum_{i=1}^d F_i$ is a  constant on  $U$. 
 
 \n These first  integrals being defined up to an additive   constant, we are only interested by their  differential $\omega_i=dF_i$, in such a way that we may still define an abelian  relation as
  a  family  $(\omega_i)_{1\leq i\leq d}$   of holomorphic 1-forms  $\omega_i$ on $U$ (maybe with singularities), which are     
 
  $(i)$      closed  (hence  locally  exact)  : $d \omega_i=0$,

   $(ii)$   verifying  $T{\cal F}_i\subset Ker \  \omega_i$   \  ($T{\cal F}_i=Ker \  \omega_i$ at any  point where  $\omega_i$ doesn't vanish), 
   
   $(iii)$ such that $\sum_{i=1}^d  \omega_i=0$.
   
   \n The  germs of abelian  relations at a    point $m$ constitute a vector space, whose  dimension is called  \emph{the  rank} of the web at this point\footnote{A. H\énaut proved that this rank doesn't  depend on $m$, as far as the web satisfies to the assumption of strong general   position   ([H2]). In case we have only weak general   position, we shall define the  \emph{rank} of the web as being the highest of the rank at a  point.}.
   
   It will be useful to give an equivalent definition in words of differential operator.
Denote by  $T{\cal F}_i \ \bigl(\subset TM\bigr)$ the vector bundle of vectors  tangent to  ${\cal F}_i$, and  $A_i\ \bigl(\subset T^*M\bigr)$ the dual vector bundle of $TM/T{\cal F}_i$ (i.e. the vector bundle of holomorphic 1-forms vanishing on  $T{\cal F}_i$).
Let   $$Tr:\oplus_{i=1}^d A_i\to  T^*M$$  be the morphism of vector bundles  (the  \emph{Trace}), defined by  $  Tr\bigl(( \omega_i)_i\bigr)= \sum_{i=1}^d  \omega_i .$  The assumption of ``at least  weak general position''  means that $Tr $ has maximal rank $n$ : its kernel  $$A:=Ker \ Tr$$ is therefore a holomorphic vector bundle of rank   $d-n$.  We define a linear differential operator of order one
$$D:J^1 A\to B, $$  where $B=(\wedge ^2T^* M)^{\oplus d}$,  by mapping any  section $s=( \omega_i)_i$ of  $A$ onto the   family $(d \omega_i)_i$ of the  differentials. Then, an   \emph{abelian relation}   may be identified with a holomorphic   section  $s$ of  $A$ such that $D(j^1s)=0$.

 
  
  
  The kernel   $R_1=Ker (D:J^1 A\to B)   $ is the vector bundle  of formal abelian relations at order one.  More generally, the space $R_h$ of \emph{formal abelian relations} at order $h$ is the kernel of the  $(h-1)^{th}$-prolongation  $D_h$ of 
  the  differential operator $D\ (=D_1)$ :
  $$R_h=Ker (D_h:J^h A\to J^{h-1}B)  .$$
 For any $h$ ($h\geq 1$), abelian relations   may still  be identified with  holomorphic   sections  $s$ of  $A$ such that $j^hs $ belong to $R_h$. 
  
  Denoting by  $\pi_h:R_h\to R_{h-1}$ the natural  projection, we shall see that the  elements of  $R_h$ which are mapped by   $\pi_h$ onto a given element $a_{h-1}$ of  $R_{h-1}$ are the  solutions of a linear  system    $\Sigma_h(a_{h-1})$ of  $c(n,h+1)$ equations with  $d$ unknowns, whose homogeneous part    doesn't depend on  $a_{h-1}$, with notation\footnote{We prefer this notation to the usual one $\begin{pmatrix} n-1+h\\h\end{pmatrix} $
 for the binomial coefficient, mainly because it suggests  explicitly the dimension of the vector space of homogeneous polynomials of degree $h$ with $n$ unknowns, and also because it needs less space.} $$c(n,h):=\frac{(n-1+h)!}{(n-1)!\ h!}.$$ Then  \emph{ordinary webs} are those for which all of these  systems have  maximal  rank  $ {inf} \bigl(d,c(n,h+1)\bigr)$.
 Denoting by  $h_0$ the integer
 such that  $$c(n,h_0-1 )< d\leq c(n,h_0),$$ 
 it is in fact sufficient that this rank be maximal for $h\leq h_0$, for being maximal for any $h$. 
 
 We proved in [CL] that the rank of an ordinary web is at most equal to the integer 
 $$\pi'(n,d):=\sum_{h=1}^{h_0-1}   \bigl(   d-c(n,h) \bigr), \ \Bigl(=(h_0-1)d-c(n+1,h_0-1)+1\Bigr).$$
  which,  for $n\geq 3$, is strictly smaller than  the  Castelnuovo's number\footnote{The Castelnuovo's  number  
  $$\pi(n,d):=\sum_{h\geq 1}\bigl(d-h(n-1)-1\bigr)^+\ , \hbox {\ where  $a^+$ denotes the number  $ sup\ (a,0)$}  $$
  is the maximal  arithmetical genus of irreducible algebraic curves of degree $d$  in the complex $n$-dimensional projective space  $\hbox{\bb P}_n$.   It is also,  after  Chern ([C]), the  maximal rank of the  $d$-webs  in codimension one,  verifying only the assumption of strong  general position (but not  necessarily ordinary for $n\geq 3$).} $\pi(n,d)$.
  
  Since the linear system $\Sigma_h(a_{h-1})$ of   $c(n,h+1)$ equations with  $d$ unknowns has   rank $c(n,h+1)$ for $h\leq h_0-2$, 
the projection $\pi_h:R_h\to R_{h-1}$ is then surjective, and $R_k\to M$ is a holomorphic vector bundle of rank 
$\sum_{h=1}^{k+1}   \bigl(   d-c(n,h) \bigr)$  for $k\leq h_0-2$. 
In particular $$R_{h_0-2}\to M\hbox{ is a holomorphic vector bundle of rank }\pi'(n,d).$$

For $h\geq h_0-1$, $\Sigma_{h}(a_{h-1})$ has   rank $d$, and has at most one solution (since it contains a cramerian sub-system), but may be no one  (since it is overdetermined). 
In general,  $R_h$ will still be  a vector bundle (denoting by $\rho_h$ its rank), but it may happen that the projection $\pi_{h+1}:R_{h+1}\to R_h$  be no more surjective, hence : $\rho_h\geq \rho_{h+1}$.
 
When $\rho_h= \rho_{h+1}$, $(h\geq h_0-1)$, the projection $\pi_{h+1}:R_{h+1}\to R_h$ is now  an isomorphism of vector bundles. The inverse isomorphism $R_{h}\buildrel\cong\over\rightarrow R_{h+1}$ composed with the natural inclusion $R_{h+1}\subset J^1R_{h}$ defines a connection $\nabla ^{h}$ on $R_h$, and   abelian relations may be identified with sections $s$ of $A=R_{0}$ such that $j^h s$ be a section of $R_h$ with vanishing covariant derivative : $\nabla ^{h} (j^h s)\equiv 0$. Hence the rank $\rho(W)$ of the web will be at most equal to the rank $\rho_h$ of $R_h$, and equal iff the curvature $K_h$ of $\nabla ^h $ vanishes. This      proves in particular the inequalities 
$${\rm sup} \Bigl(0,(h+1)d-c(n+1, h+1)+1\Bigr)\leq \rho_{h}\leq (h_0-1)d-c(n+1, h_0-1)+1=\pi'(n,d)$$for $h\geq h_0-2$. We shall see how to compute effectively $\rho_{h}$ and $K_{h}$.

  When $d=c(n,h_0)$  (then, $d$ is said to be \emph{calibrated}), $\rho_{h_0-2}=\rho_{h_0-1}$, and  $\nabla^{h_0-2}$  is the connection defined in [CL], computed with Maple in [DL], generalizing to any $n$ the connection previously defined for  planar webs ($n=2$) by  H\énaut ([H1]), and independently by Pirio ([Pi]) who related  the corresponding curvature to invariants defined previously by Pantazi ([Pa]) (the  curvature of which being the Blaschke-Dubourdieu curvature ([BB])  for $n=2,\ d=3$). 
In this case of  planar webs, Ripoll ([R]) computed   the rank of the web by another method (corank of some matrix deduced from   $\nabla^{h_0-2}$, from the curvature  $K^{h_0-2}$ and    its derivatives).

\n {\bf Remark about non-ordinary webs :}  If the web is not ordinary, its rank may be bigger than $\pi'(n,d)$, and we  may not affirm anymore that the sequence of the $\rho_h$'s increases for $h\leq h_0-2$ and decreases for $h$ bigger.   
 However, if -by chance- we  can find some $h$ such that  $\pi_{h+1}:R_{h+1}\to R_h$ is  an isomorphism of vector bundles of same rank $\rho_h=\rho_{h+1}$, then we still can define the connection $\nabla^h$, and it is still true that the vanishing of the curvature $K^h$ implies the equality $\rho(W)=\rho_h$.
 
\n {\bf Notation :} Denote by

 \n $i$ an index from  1 to $d$,
 
 \n $\lambda,\mu,...$ an index from  1 to $n$,

 \n  $L=(\ell_1,\ell_2,\cdots,\ell_n)$ a  multi-index  $\ell_\lambda\geq 0$ of $n$ integers, and $|L|:=\sum_\lambda \ell_\lambda$  its  \emph{degree}.  
  
  \n If $L=(\ell_1,\ell_2,\cdots,\ell_n))$, and $L'=(\ell'_1,\ell'_2,\cdots,\ell'_n))$, $L+L'$ (resp. $L-L'$)  denotes $(\ell_1+\ell'_1,\cdots,\ell_n+\ell'_n)$ (resp $(\ell_1-\ell'_1,\cdots,\ell_n-\ell'_n)$).

 \n In particular $ 1_\lambda$ denotes the  multi-index obtained with  1 at the  place $\lambda$ and  0 elsewhere.
    
  \n Relatively to local  coordinates   $x =(x_1,\cdots,x_n)$    in   $M$,  we shall  denote by $\partial_\lambda a$ or $a'_\lambda$ the  partial   derivative 
  $\frac{\partial a}{ \partial x_\lambda } $ of a holomorphic function   $a$ or of a    matrix with holomorphic coefficients.  \n More  generally, $\partial_L a$ or $a'_L$ denotes   the  partial   derivative $\frac{\partial^{|L|}a}{(\partial x_1)^{\ell_1}\cdots(\partial x_n)^{\ell_n}} $ of order $|L|$.

 
 \section{Computation of $R_h$}
  
  We assume each foliation  ${\cal F}_i $ defined by a first integral $u_i$ without  singularity.
The data of  another first integral
$F_i=G_i(u_i)$ up to an   additive constant  is equivalent to the data of the   derivative $g_i=(G_i)'$.  Each  vector bundle  $A_i$ being  now  trivialized by $du_i$, we set   $\omega_i=g_i(u_i)\ du_i$ (such a 1-form is  automatically closed). The data of an 
 abelian  relation  is now equivalent to the data of a  family   $(g_i)$  of holomorphic   functions of one  variable  ($1\leq i\leq d$) such that  $\sum_i g_i(u_i)\ du_i\equiv 0$, 
or equivalently :  $$(E_\lambda)\hskip 1cm  \sum_i (u_i)'_\lambda\ g_i(u_i)  \equiv 0\hbox {\hskip 1cm  for any $ \lambda$},$$
which can still be written $$<{P}_1\ ,\ f>\equiv 0,$$ where $P_1:=\frac{D(u_1, \cdots,u_d)}{D(x_1,\cdots,x_n)}$ denotes the jacobian matrix and $f$ the $d$-vector $\bigl(g_1\scirc u_1,\cdots,g_d\scirc u_d\bigr)$. 
 [The  functions $u_i$  are given, and the  functions $g_i$  unknown].

\n{\bf Coefficients  $C^{h}_{ L}( u )  $ and matrices ${M}_j^{(h)}$}

    For any  $h\geq 0$, $g_i^{(h)}$ will denote the  $h^{th}$  derivative of  $g_i$ (with the  convention $g_i^{(0)}:=g_i)$ ;
  
  We set :
  
   $f_i:=g_i\scirc u_i$ and   $f_i^{(h)}:=g_i^{(h)} \scirc u_i$, 
    
$ f:=$ $d$-vector $(f_1,f_2,\cdots,f_d)$, and   $ f^{(h)}:=$ $d$-vector $(f_1^{(h)},f_2^{(h)},\cdots,f_d^{(h)})$.
  
 For any integer  $k$ ($k\geq 0$), a $k$-jet of abelian  relation at a  point $m$ of $M$ is  defined by the  family $$\Bigl(f_i^{(h) }(m)=(g_i^{(h)}\scirc u_i)(m)\Bigr)_{i,h},\ \ (0\leq h\leq k,\ 1\leq i\leq d).$$
 \n 
The     partial  derivatives of the  relations $(E_\lambda)$ will make us able  to compute locally  $R_h$. In fact, 
the  functions $C_{i,{L}}^{h}$ will be defined by iteration on $|L|$ so that 
$$\Bigl(f_i.(u_i)'_\lambda\Bigr)'_L\equiv \sum_{h=0}^{|L|} C_{i,{L+1_\lambda}}^{h}\ .\ f_i^{(h)}$$as far as  $(f_i\ du_i)_i$ be an abelian  relation.

\n {\bf Lemma 2-1 :} {\it For any     holomorphic function  $u$ of  $n$ variables, and any holomorphic  function    $g$ of one  variable, 

$(i)$ The  derivatives $\Bigl((g\scirc u)\ u'_\lambda\Bigr)'_L$  are linear combinations  
$$\Bigl((g\scirc u)\ u'_\lambda\Bigr)'_L=\sum_{h=0}^{|L|} C_{L+1_\lambda}^h(u)\ .\ (g^{(h)}\scirc u)  $$ of the   successive derivatives $g^{(h)}$ of  $g$ $($we set : $g^{(0)}=g)$, whose  coefficients  $C_{L' }^h(u) =C_{L+1_\lambda}^h(u)$ depend only on  $u$  and on the multi-index $L'=L+1_\lambda$, and not on its  decomposition   under the shape  $L+1_\lambda$.

$(ii)$ They can be computed by iteration on    $|L|$ using  the formula

 $$\begin{matrix}
  C^0_{1_\lambda} (u)  & =&\ u'_\lambda\ ,&&\\
    &&&&&  \\
 
  C^0_{  L+1_\mu}(u)&=& \partial_\mu C^0_{  L}(u)\ , & &\\
 
  &&&&&  \\
 
  C^h_{  L+1_\mu}(u)&=& \partial_\mu C^h_{  L}(u)\  +&\ C^{h-1}_{ L} (u)  \ .\  u'_\mu  &\hbox{ for }  1\leq  h\leq |L|-1\ ,\\

  &&&&&  \\

C^{|L|}_{L+1_\mu}(u)&=& C^{|L|-1 }_{ L} (u)  \ .\  u'_\mu \ .&& & \\

\end{matrix}  $$ }

\n  This  is due to the fact that the 1-form $d\bigl(G(u)\bigr)$ is closed, $G$ denoting a  primitive of  $g$.

\n For a web  locally defined by the functions    $u_i$, we set :  $$C^h_{i,  L}=C^h_{ L}(u_i) .$$

\n We check in particular 
 
 $C_{i,{L}}^{0}=(u_i)'_L$, 
  
  and 
 $ C_{i,{L}}^{|L|-1}=\prod_{\lambda=1}^n \bigl((u_i)'_\lambda\bigr)^{\ell_\lambda} \hbox{ for  $L=(\ell_1,\ell_2,\cdots,\ell_n)$}$.

 \n We set :

  $\Theta ^r$ denotes the trivial holomorphic bundle of rank $r$,
  
  $  \beta_k:=c(n+1, k)-1 \ \ \bigl( =\sum_{h=1}^kc(n,h) \bigr),$
  
   ${M}_j^{(h)}$ denotes the matrix  $(\!(C_{i,{L}}^{h})\!)_{(i,|L|=j)}$ of size  $c(n,j)\times d$,  $(1\leq j$, $0\leq h)$,   \hb \indent (with ${M}_j^{(h)}=0$ for  $h\geq j$), 
 
  $P_j:={M}_j^{(j-1)}$.

 ${\cal M}_k$ denotes the matrix of size $\beta_k\times kd$ built with the  blocks ${M}_j^{(h)}$ for  
 $1\leq j\leq k$ and   $0\leq h\leq k-1$,\hb 
 \indent   where the  block ${M}_j^{(h+1)}$ is on the right of  ${M}_j^{(h)}$, and  ${M}_{j+1}^{(h)}$ below,

  $Q_{k+1}$ denotes the sub matrix of size $c(n,k+1)\times kd$ in ${\cal M}_{k+1}$ built with the  blocks ${M}_{k+1}^{(h)}$ for     $0\leq h\leq k-1$ :


 $$\begin{matrix}
 
 {}{\cal M}_k  &  =  & 
 
 \begin{pmatrix}{M}_1^{(0)}=P_1&0&0&....&....&0&0\\ 
 {M}_2^{(0)}& {M}_2^{(1)}=P_2&0&....&....&0&0&\\
....& .... &....&....&....&0&0&\\ 
  {M}_{k-1}^{(0)}& {M}_{k-1}^{(1)} &{M}_{k-1}^{(2)}&....&....&{M}_{k-1}^{(k-2)}=P_{k-1}&0\\ 
{M}_k^{(0)}& {M}_k^{(1)} &{M}_k^{(2)}&{M}_4^{(3)}&....&{M}_k^{(k-2)}&{M}_k^{(k-1)}=P_k\\ 
 & \end{pmatrix}
 \\
 &&\\
Q_{k+1}  &  =  &  

 \begin{pmatrix} \hskip .3cm
 
{M}_{k+1}^{(0)}\hskip .5cm& {M}_{k+1}^{(1)}\hskip .3cm &{M}_{k+1}^{(2)}\hskip .2cm&{M}_{k+1}^{(3)}\hskip .1cm&....\hskip .5cm&{M}_{k+1}^{(k-2)}\hskip 1.2cm&{M}_{k+1}^{(k-1)}\hskip .8cm\\ 
 
& \end{pmatrix}
 
 \end{matrix}$$


 \n {\bf Theorem 2-2 :} 
 
 {\it 
 
 $(i)$ Locally, $R_k$ is the kernel of ${}{\cal M}_{k+1}$ $($included  into the trivial bundle $\Theta^{(k+1)d})$. Hence, when the matrix ${}{\cal M}_{k+1}$ has constant rank $($always true for $k\leq h_0-2)$, $R_k\to M$ is a holomorphic vector bundle of rank $$\rho_k=(k+1)d-rank({}{\cal M}_{k+1}).$$
 
 $(ii)$ If an element $a_{k-1}\in R_{k-1}$  (imbedded into $\Theta ^{k d}$) is defined by the family $f$ of  $d$-vectors    $  (f^{(h)})_h   $, $(0\leq h\leq k-1)$, the elements of   $  R_{k }$ which project onto   $a_{k-1}$ are those whose the  last component  $f^{(k)} $ is  solution  of the linear system $\Sigma_h(a_{h-1})$ :
 $$ \begin{matrix}  <P_{k+1}\ , f^{(k)}> &=&-  <Q_{k+1}\ ,\ f> \\&&\\
 &=&-\sum_{h=0}^{k-1}<{M}_{k+1}^{(h)}\ ,\ f^{(h)}> \end{matrix}.$$}
 
\n {\it Proof :}
 A $k$-jet of abelian relation  at a point $m\in M$ is then     represented by its  components  $ j_m^h(u_i\scirc g_i )$ in  $J^h A_i$, and each of them is completely defined by the  family of the numbers     $$ \Bigl( f_i^{(h)}=(g_i^{(h)}\scirc u_i)(m)\Bigr)_{0\leq h\leq k} .$$  Thus  a family\footnote{Be careful not to confuse $R_h$ with the set $J^h(Ab)$ of the $h$-jets of the true abelian relations (which may be smaller).}  of numbers   $\bigl( w_i^{(h)}\bigr)_{i,h}$ belongs to  $R_k$ if it  satisfies to any of the equations 
 $$(E_L)\hskip 1cm  \sum_{i=1}^d\sum_{h=0}^{|L|-1} C_{i L}^h\ .\ w_i^{(h)}  =0.                       $$

   \rightline{QED}
  
 \n {\bf Estimation of the ranks}
 
  \n The assumption for the web to be ordinary means  that 
 the  matrices  $P_j:={M }_j^{(j-1)}$ have all  maximal rank, that is $c(n,j)$ for $j\leq h_0-1$, and $d$ for $j\geq h_0$.

 \n {\bf Lemma 2-3}
 
  \n {\it If $P_j$ has maximal rank $c(n,j)$ for $1\leq j\leq h_0$, it has maximal rank $d$ for any $ j\geq h_0$.}
  
  \n {\it Proof :} The meaning of the lemma not depending on the local coordinates, we may assume that all foliations are transversal to the $x_n$-axis near a point ; therefore all derivatives $(u_i)'_n$ are not zero. The formula 
  $$\bigl(P_{j+1}\bigr)_{i,L+1_n}=(u_i)'_n.\bigl(P_{j}\bigr)_{i,L}$$
  proves that the rank of $P_{j+1}$ is at least equal to that of $P_j$, thus is equal if  $j$ is  big  enough for this rank to be stationary equal to $d$. 
  
  \rightline{QED}
  
  Let $\rho_k$ be the rank of $R_k$.
  
  \n {\bf Theorem 2-4 ([CL])}
 
\n {\it  For  $k\leq h_0-2$,   $R_k\to M$ is a holomorphic vector bundle of rank  
 $$\rho_k= (k+1)d-\beta_{k+1} .$$
 In particular, $\rho_{h_0-2}=\pi'(n,d)$.}
 
   \n {\it Proof :} In fact, the matrices ${}{\cal M}_{h}$,    of size $\beta_h\times hd$, are triangular by blocks, and the diagonal blocks are the $P_j$'s. 
   Since the rank of $P_j$ is $c(n,j)$ for $j\leq h_0-1$, ${\cal M}_{k+1}$ has maximal rank $\beta_{k+1}=\sum_{h=1}^{k+1} c(n,h)$ in this range. Thus $R_k \ (={\rm Ker} \ {}{\cal M}_{k+1})$ has there rank $(k+1)d-\beta_{k+1} $.
   
     \rightline{QED}

The sequence $(hd-\beta_h)_h $ becomes decreasing  for $h\geq h_0$
Then, for  $h\geq h_0$, it may be no more true that ${}{\cal M}_{h+1}$ has maximal rank $\beta_{h+1}$, so that the rank $\rho_h$ of $R_h$ may be now bigger than $(h+1)d-\beta_{h+1} $ (but remains at most equal to $\pi'(n,d)= (h_0-1)\ d-\beta_{h_0-1} $). Thus we get :

\n {\bf Theorem 2-5 }

{\it Assuming $\pi_{h+1}$ to have a constant rank for $h\geq h_0-2$, the sequence $(\rho_h)_{h\geq h_0-2}$ is decreasing from $\pi'(n,d)$ to the rank $\rho(W)$ of the web, and satisfies to the inequalities 
$$(h+2)d-\beta_{h+2} \leq \rho_{h+1}\leq \rho_{h}\leq \pi'(n,d).$$}


 
  }

\section{The connections $\nabla^h$}

In this section, we assume :

$h\geq h_0-2$,  
 
 and  $\rho_h=\rho_{h+1}$, \hskip 1cm  ($\pi_{h+1}:R_{h+1}\buildrel\cong\over \rightarrow R_h$ being then   an isomorphism of vector bundles). 
 
\n  If $\rho_h=0$, then $\rho(W)=0$. If $\rho_h>0$, we shall   define a connection $\nabla^h$ on   $R_h$, whose curvature vanishes iff $\rho_h=\rho(W)$.   
 
\n We recall that 
 $ R_{h+1}$ is the intersection of  $J^1(R_{h})$ and $ J^{h+1}R_0$ into $J^1(J^h R_0)$  :  $$ R_{h+1}=J^1 R_{h} \cap J^{h+1}R_0.$$ 
 
  \n Denote by 
  
  $\epsilon_h: R_{h+1} \hookrightarrow J^1 R_{h}$ the natural inclusion, 
  
  and by $v_h: R_h\to R_{h+1}$ the   isomorphism inverse\footnote{It can be explicitely computed by mean of the generalized inverse $IP:=(P^*.P)^{-1}.P^*$ of $P$\ , where $P:=P_{h+2}$ and $P^*$ means the transposed matrix of $P$.}  of the projection $\pi_{h+1}$.
  
\n The composed map  $ \xi_h:=\epsilon_h\scirc v_h$ is a splitting of the exact sequence
$$0\to T^*(V)\otimes R_h\to J^1 R_h\buildrel{\buildrel { \xi_h} \over \longleftarrow}\over\longrightarrow R_h\to 0  $$ 
and   defines consequently a holomorphic connection on  $R_h$, whose covariant derivative is  :
$$\nabla^h  \sigma =j^1\sigma-<\xi_h,\sigma>.$$
Since the abelian  relations may be identified by the map  $s\to j^{h+1}s$ to the sections $s$ of  $R_0\ (=A)$ such that  $j^{h+1}s$ belong to   $R_{h+1}$,  
and since $\xi_h$ factorizes through $R_{h+1}$, the following assertions are equivalent :

$(i)$ \ $s$ is an abelian relation,

$(ii)$ \ $\nabla ^h (j^{h}s)\equiv 0$.

\n Since the framework is holomorphic, $\rho_\infty=\rho(W)$, and  we get therefore the 

\n {\bf Theorem 3-1 } {\it 

$(i)$ \ \ A section  $s$ of  $A\ (=R_0)$ is an abelian  relation iff $j^{h}s$ is a section of $R_h$ and $\nabla^h (j^{h}s)\equiv 0$.

$(ii)$ \ The  rank $\rho(W)$ of the web is at most equal to the rank  $\rho_{h}$ of the bundle  $R_h$.  

$(iii)$  There exists an integer $h_1$ such that 

-  or $\rho_{h_1}=0$ and then $\rho(W)=0$, 

- or  $\rho_{h_1}=\rho_{h_1+1} \ (\neq 0)$, the   curvature $K^{h_1}$ vanishes , and then  $\rho(W)=\rho_{h_1}.$ }
\n {\bf Remark : } If the web is  not-ordinary, we still may define the connection  $\nabla ^h$, as soon   as we can find some $h$ for which the    projection $R_{h+1}\to R_h$ is an isomorphism of vector bundles, whatever be $h$ (no more necessarily $\geq h_0-2$).  And it remains true that the vanishing of its curvature $K^h$ implies $\rho(W)=\rho_h$. 
\section{Algorithm}
Theoretically, the following algorithm always works for any ordinary web. But it may need  a long  time of computer.
Practically, in some cases,   considerations specific to each example may be used for shortening the process, some of them being sketched in the remark at the end of the section.

- explicit  $P_1,\cdots, P_{h_0}$  ; 

- check : $${\rm Rank}\ (P_j)=c(n,j)   \hbox{ for any   $j$   ($1\leq j\leq h_0-1$)}\hbox {, and } {\rm Rank}\ (P_{h_0})=d\ ;$$     \indent  if this condition is not realized,  $W$ is not ordinary ; 

- STOP ;

- else,  compute ${\cal M }_{h_0}$ ;   

\n  {\bf Loop ${\cal L}$(h) from $h=h_0-2$ :}

- compute ${\cal M }_{h+2}$ (and ${\cal M }_{h+1}$ sub-matrix of  ${\cal M }_{h+2}$) ;

- compute $\rho_{h}=(h+1)\ d-{\rm Rank}\ ({\cal M }_{h+1})$ and $\rho_{h+1}=(h+2)\ d-{\rm Rank}\ ({\cal M }_{h+2})$ ;

- if $\rho_{h }>\rho_{h+1}$, go to ${\cal L}$(h+1) ;

- else (i.e. when $\rho_{h }=\rho_{h+1}$), compute $\nabla ^h$ and $K^h$ ;

- if $K^h\neq 0$, go to ${\cal L}$(h+1) ;

- else  (i.e. when $K^h= 0$),   $$\rho(W)=\rho_h.$$

-  STOP .

\n Thus,  we have an effective procedure to compute the rank, even  when it is not maximal, without having to exhibit explicit abelian relations.

\n {\bf Remarks  :} 

1-  When $\rho_h=\rho_{h+1}$,   it is often useful to check immediately if there would not be  some $k$, $k>h$, such that $\rho_k>\rho_{k+1}$. In this case, we know a priori that $K^h$ doesn't vanish, without to have to compute it. 

2- There are usually two ways for computing $\rho_h$ : the first one, used in the algorithm above, consists in computing the kernel of the matrix  ${\cal M }_{h+1}$ of size 
 $\bigl(c(n+1,h+1)-1\bigr)\times (h+1)d$ : 
 $$\rho_h=(h+1)\ d-{  Rank} ({\cal M }_{h+1} ) .$$This  size 
  increases more rapidly with $h$ than the size $c(n, h+1)\times d$ of the matrix $P_{h+1}$ of the linear system 
$\Sigma_h(a)$ giving the elements of $R_h$ above a given element $a\in R_{h-1}$
 (essentially because the process uses   the knowledge of $R_{h-1}$ that we got previously, which is not true for the first process). Thus, for $h$ big enough, knowing already $\rho_{h-1}$ and a trivialization $(\epsilon_s)_s$ of $R_{h-1}$, the following process may need a   shorter time of computer than the previous one, despite of the fact that there are more operations to be done : 

- choose a $d\times d$ invertible sub-matrix  $P^0_{h+1}$ of $P_{h+1}$,

- solve the corresponding cramerian sub-system of $\Sigma_h(a)$, 

- for each line $\ell$ among the $c(n,h+1)-d$ deleted for getting $P^0_{h+1}$ from $P_{h+1}$, and for each $\epsilon_s$ belonging to the trivialization  of $R_{h-1}$, build the characteristic determinant $\Delta(s, \ell)$ whose vanishing asserts the compatibility  of the new equation $\ell$ with the cramerian sub-system,

- then the kernel of the matrix $\Delta^h:=(\!(\Delta(s, \ell))\!)$ of size $\bigl(c(n,h+1)-d\bigr)\times\rho_{h-1}$  defines the projection of $R_{h}$ onto $R_{h-1}$, which is an isomorphism, and $$\rho_h=\rho_{h-1}-{  Rank} (\Delta^h).$$ 

 \vskip -2 cm

\section{Examples}

The process described in the algorithm above works for any $(n,d,h)$. However,  most of our examples are relative to low values of these integers : in fact,  the size of the involved matrices becomes  very rapidly huge, and would need in practice more powerful computers  than our small portable.

\subsection{Case $n=2,\ d=3$ :}
 There is no hope to refine the classification of the non-hexagonal  planar 3-webs by the order of the step from which the sequence of the $\rho_h$'s vanishes, In fact, 
we can prove easily that
 the sequence of the $\rho_h$'s becomes immediately stationary after the first step, and     there are only two possibilities :

- sequence $(1,1,\cdots,1=\rho_\infty)$ if the Blaschke-Dubourdieu curvature $K^0$ vanishes  $($hexagonal case$)$,

-  sequence  $(1,0,\cdots,0=\rho_\infty)$ if $K^0\neq 0$ .

\subsection{Example $n=2,$ $ d=4\ $$\bigl(\pi'(2,4)=3\bigr)$ :  }
 
 We recall that all  planar  webs   are 
 ordinary,
and calibrated  with $h_0=d-1$. Moreover $\pi'(2,d)$ is then equal to $\pi(2,d)=\frac{(d-1)d-2)}{2}$.

\n For  the  planar  4-web   $$(x,y,x+y+xy,x-y+x^5), $$ 
we have an obvious abelian relation $f\scirc u_1-u_2-u_4\equiv 0$, with $f(x):=x+x^5$. Thus, we know already $$1\leq \rho(W)\leq 3.$$
Computing $\rho_k$, we get $$\rho_1=\rho_2=3 >\rho_3=\rho_4=2.$$Since $\rho_3<\rho_2$, we are sure that the curvature $K^1$ doesn't vanish,  without to have  to compute it. We get $K^3=0$. Thus 
the sequence of the $\rho_i$'s  is necessarily  stationary equal to 2 from $\rho_3$ :
$$\rho(W)=2.  $$
We are sure that there is another   abelian relation independant on  the obvious one, without to have to exhibit it.

\subsection{Example  $n=2,$ $ d=8\ $ $\bigl(\pi'(2,8)=21\bigr)$ : }

Let $W$ be the planar $8$-web 
$$x,y,x+y,x-y, xy,x^2+y^2,x^2-y^2,x^4+y^4.$$
We observed in [DL] that its curvature $K^{5}$ did not vanish, but that its  connection form $\omega^5$ relative to some ``adapted'' trivialization (matrix of size $21\times 21$, whose coefficients are scalar 1-forms) had only zero's in the 19$^{th}$ first columns. Thus, we deduced that the rank of $W$ was at least 19, and at most 20 :
in fact, $\rho_5=21,\ \rho_6=20$ and \ $\rho_7=19$. Thus 
$$\rho(W)=19,$$
while the $7$-sub-web defined by deleting $x^4+y^4$ has maximal rank 15 (see [Pi]).

\subsection {Case  $d=n+1$, $n>2$\ $\bigl(\pi'(n,n+1)=1\bigr)$ : }

Denoting by $(x_1,\dots ,x_n)$ local coordinates, we consider the $(n+1)$-web $W$ defined by the functions $(x_1,\dots ,x_n,F(x_1,\cdots ,x_n)).$

\n For a convenient  order of the multi-indices $L$, the matrix ${\cal{M}}_2$ has the shape
$$\begin{pmatrix}     I_n  & F^{(1)} &  0  & 0  \\ 
                        0   &F^{(2)}  & I_n & F^2 \\ 
                        0 &  G^{(2)}  & 0   & G^2    \end{pmatrix},$$ 
where 

$I_n$ is the identity $n\times n$-matrix  , 

$F^{(1)}$ is the column of the $(F'_i)_{1\leq i\leq n}$, (with $n$ rows), 

 $F^{(2)}$ is the column of the $(F''_{ii})_{1\leq i\leq n}$, (with $n$ rows),

 $F^2$ is the column of the $((F'_{i})^2)_{1\leq i\leq n}$, (with $n$ rows),

 $G^{(2)}$  is  the column  which has   coefficients $F''_{ij}$ ($i\neq j $, (with $c(n,2)-n$ rows),
 
 and   $G^2$ is  the column  which has   coefficients $F'_i F'_j$  ($i\neq j $, (with $c(n,2)-n$ rows,
 same order of the \hb \indent pairs $(i,j)$ as in $G^{(2)}$).
 
\n The sub-matrix ${\cal{M}}_1$ has always   rank $n$, hence $\rho_0=1$,  while ${\cal{M}}_2$ has generally rank $2n+2$ ; hence, in general $\rho_1\ \bigl(=\rho(W)\bigr)=0$, and  there is no abelian relation.

\n The exceptional case   (Rank$({\cal{M}}_2)=2n+1$, and $\rho_1=1$) happens iff $G^{(2)}$ and   $G^2$ are collinear.
 This means  the set of relations
$${{F''_{ij}}\over{F'_iF'_j}}\equiv {{F''_{rs}}\over{F'_rF'_s}}$$
for any $i,$ $j,$ $r$ and $s$ with $j\neq i$ and $s\neq r.$

\n We shall now study this case 
by mean of the connection $\nabla ^0$.  
Then, a trivialization of  $R_0={\rm Ker}\ {\cal M}_1$ is given  by the $(n+1)$-vector  $$f^{(0)}=(-F'_1,-F'_2,\cdots,-F'_n,1)$$ and
a trivialization of  $R_1={\rm Ker}\ {\cal M}_2$ is given  by some  $(n+1)$-vector  $$f^{(1)}=(X_1,X_2,\cdots X_{n+1})$$ 
 satisfying in particular to the identities
$$X_{n+1}=-{{F''_{ij}}\over{F'_iF'_j}} \hbox{ \ \  whatever be } i,j,\ (i\neq j).$$

\n Denoting by $\Delta_i$ the $(n+1)\times (n+1)$ diagonal matrix built on the $(n+1)$-vector $$(0,\cdots,0,1,0,\cdots,0,F'_i),$$ with $1$ as $i^{th}$ component, $F'_i$ as  $(n+1)^{th}$ component and 0 elsewhere,  the connexion $\nabla^0$ on $R_0$ is then defined by
$$\nabla^0_i f^{(0)}=\partial_i f^{(0)}-<\Delta_i, f^{(1)}>,$$
where $\partial_i$ (resp. $\nabla^0_i$) means the partial derivative (resp. the covariant derivative) with respect to 
$\frac{\partial}{\partial x_i}$. 
Thus, we get :
$$\nabla^0_if^{(0)}=-F'_i X_{n+1} f^{(0)}.$$
The curvature has then components
$$K^0_{ij}=\partial_i(F'_j X_{n+1})-\partial_j(F'_i X_{n+1}).$$
Fix a pair  $(i,j)$ and choose an index $k$  different from $i$ and $j.$ We get :
$$F'_iX_{n+1} =-F''_{ik}/F'_k,$$
hence 
$$\begin{matrix}\partial_j(F'_iX_{n+1}))&=&-\partial_j(F''_{ik}/F'_k)\\&=&-\partial_j(\partial_i\ln (F'_k))\end{matrix}.$$
This gives  $K_{ij}^0=0$. So, if we set 
$$L_{ij}=ln(F'_{i}/F'_{i}), $$  we have proved   the following proposition.
\n {\bf Proposition 5-2 :} 

\n {\it The web $W$ has a non-trivial abelian relation if and only iff $$\ (L_{ij})'_{k}=0,$$
for any triple $i,j,k$ of indices,  each one   different to each other.}

\n Notice that, when $W$ is in strong general position, the existence of an abelian relation is equivalent to the fact that we can choose new coordinates $(\overline{x}_i)_i$ such that 
$$F(x_1,\dots , x_n)\equiv \overline{x}_1+\cdots +\overline{x}_n.$$
Thus,  the existence of an abelian relation is equivalent for the web to be ``parallelisable''.




\subsection{An example  $n=3,\ d=5$\ $\bigl(\pi'(3,5)=2\bigr)$ :}
Denoting by $(x,y,z)$ local coordinates, and  defining  the web  by the  
functions $\bigl(x,y,z, x+y+z, F(x,y,z)\bigr)$, assume that the function $F$ depends only on $ x+y$ and $z$ :
$$ F(x,y,z)\equiv g(x+y,z)\hbox{ for some function } g.$$We set 
$u:= x+y,$

$p:=g'_u,\ q:=g'_z,$ 

$r:=g''_{u^2},\ s:=g''_{u z}, \ t:=g''_{z^2},$ 

$a:=g'''_{u^3},\ b:=g'''_{u^2z},\ c:=g'''_{u z^2},\ e:=g'''_{z^3},\ $.

 \n We consider $ {}{\cal M}_1$,  $ {}{\cal M}_2$, $Q_3$ and $P_3$ as sub-matrices of the matrix $ {}{\cal M}_3$  described below for a convenient order of the multi-indices $L$ :
$$
 {}{\cal M}_3  =
 \begin{pmatrix}     \begin{pmatrix}    1&0&0&1&p\\0&1&0&1&p\\ 0&0&1&1&q\end{pmatrix}                  &0&0\\&&\\ 
 \begin{pmatrix}    0&0&0&0&r\\ 0&0&0&0&r\\  0&0&0&0&t\\ 0&0&0&0&s\\ 0&0&0&0&s\\ 0&0&0&0&r\end{pmatrix}  &\hskip -.4cm\begin{pmatrix}    1&0&0&1&\hskip .5cm p^2\\ 0&1&0&1&\hskip .5cm p^2\\  0&0&1&1&\hskip .5cm q^2\\ 0&0&0&1&\hskip .5cm pq\\0&0&0&1&\hskip .5cm pq\\ 0&0&0&1&\hskip .5cm p^2\end{pmatrix}
 &0\\&&\\
 
 \begin{pmatrix}    0&0&0&0&a\\ 0&0&0&0&a\\  0&0&0&0&e\\ 0&0&0&0&b\\ 0&0&0&0&b\\ 0&0&0&0&b\\ 0&0&0&0&c\\ 0&0&0&0&c\\ 0&0&0&0&a\\ 0&0&0&0&a\end{pmatrix}&
  \begin{pmatrix}    0&0&0&0&3pr\\ 0&0&0&0&3pr\\  0&0&0&0&3qt\\ 0&0&0&0&2ps+rq\\ 0&0&0&0&2ps+rq\\ 0&0&0&0&2ps+rq\\ 0&0&0&0&2qs+pt\\ 0&0&0&0&2qs+pt\\ 0&0&0&0&3pr\\ 0&0&0&0&3pr\end{pmatrix}&
 \begin{pmatrix}    1&0&0&1&p^3\\ 0&1&0&1&p^3\\  0&0&1&1&q^3\\ 0&0&0&1&p^2q\\0&0&0&1&p^2q\\ 0&0&0&1&p^2q\\ 0&0&0&1&pq^2\\ 0&0&0&1&pq^2\\ 0&0&0&1&p^2q\\ 0&0&0&1&p^2q\end{pmatrix}&\\ 
\\ 
 & \end{pmatrix} $$
 We can check that $ {}{\cal M}_1$, $ {}{\cal M}_2$  above have respectively  rank 3, 8   ; thus  $$\rho_0=2\ (=5-3),\hbox { and } \rho_1=2\ (=10-8).$$  
 In general $ {}{\cal M}_3$ has rank $14$ and $\rho_2=1\ (=15-14)$. 
 But it may happen that $ {}{\cal M}_3$ has rank $13$ and $\rho_2=2$ for exceptional $g$'s. 
This can  be seen 
 by   computing   the curvature $K_0$.

\n  A basis for  $R_0={\rm Ker}\  {}{\cal M}_1$ is
 $$f_1=\Bigl(-1,-1,-1,1,0\Bigr) \ \ ,\ \  f_2=\Bigl(-p,-p,-q,0,1\Bigr).$$
The lines 7 and 8 of ${}{\cal M}_2$ being the same, we may ignore the line 8  in the computation of \hb 
 $R_1={\rm Ker}\  {}{\cal M}_2$. We assume $p\neq q$, in such a way that the   sub-matrix $P^0_2$ of $P_2$ that we get in forgetting its line 5   is invertible. Thus, $R_1$ has rank $\rho_1=2$, and we can   lift $f_1$, and  $f_2$ in $R_1$, defining $f_1^{(1)}=-<(P^0_2)^{-1}.\ {M}_2^{(0)}\ ,\ f_1>$, and  $f_2^{(1)}=-<(P^0_2)^{-1}.\ {M}_2^{(0)}\ ,\ f_2>$. We get :  
 $$f_1^{(1)}=\Bigl(0,0,0,0,0\Bigr) \ \ ,\ \  f_2^{(1)}=\Bigl(0,0,Z,T,U\Bigr) ,  $$
 where $Z$, $T$ and $U$ are solution of the cramerian linear system 
 $$\begin{matrix}
 &&T&+&p^2U&+&r&=&0\\&&T&+&pqU&+&s&=&0\\Z&+&T&+&q^2U&+&t&=&0
 \end{matrix}$$
Denoting respectively by  $\Delta_x$, $\Delta_y$, and $\Delta_x$, the $5\times 5$ diagonal matrices built with $(1,0,0,1,p)$, $(0,1,0,1,p)$,  and $(0,0,1,1,q)$, 
the connection $\nabla^0$ on $R_0$ is then given  by the formulae : 

 \hskip .4cm$ \nabla^0 f_1 \equiv 0 , $ 
 $$ \nabla^0_x f_2\equiv  \frac{\partial}{\partial x} f_2-<\Delta_x,\ f_2^{(1)}>,  \nabla^0_y f_2\equiv  \frac{\partial}{\partial y} f_2-<\Delta_y,f_2^{(1)}>,\  \nabla^0_z f_2\equiv  \frac{\partial}{\partial z} f_2-<\Delta_z,f_2^{(1)}>,$$
  where $\nabla^0_x $, $\nabla^0_y $ and $\nabla^0_z $ denote the covariant derivative with respect to $\frac{\partial}{\partial x}$, $\frac{\partial}{\partial y}$, and  $\frac{\partial}{\partial z}$.
  The connection form relative to $(f_1,f_2)$ is then
  $$\omega_0=\begin{pmatrix}0&-T(du+dz)\\&&\\0&-U\bigl(p\ du+q\ dz\bigr)\end{pmatrix}, $$
 and the curvature writes 
 $$K^0 =\begin{pmatrix}0&T'_z-T'_u+(q-p)TU\\&&\\0& pU'_z-qU'_u\end{pmatrix} (dx+dy)\wedge dz. $$ 
 If this curvature vanishes
 (according to $g$), $\rho(W)=2$. Otherwise, $\rho(W)=1$. (The rank may not be zero, due to the obvious non-trivial abelian relation  $(x)+(y)+(z)-(x+y+z)\equiv 0$). 
 For example, if $g(u,z)=u^2+2\lambda uz+\mu z^2$, ($\lambda,\mu\in ${\bb C}), we can affirm that there is no other independant abelian relation if $\lambda\neq 1$. If $ \lambda=1$, we have a vanishing  curvature, corresponding to the second abelian relation
 $u_5\equiv (u_4)^2+(\mu-1)(u_3)^2$.

\subsection{An example  $n=3,\ d=11$\ $\bigl(\pi'(3,11)=14\bigr)$  :}

Let $W$ be the $11$-web (quasi-parallel : all $u_i$'s except one are affine  functions) :
$$x,\ y,\ z,\   x+y+z,\ x+2y+z,\ x+3y+z,\ x+y+5z,\ x+y+7z,\ x+11y+z,\ 19x+y+z,\  x+yz.$$
We get $\rho_2=14>\rho_3=\rho_4=13$, and $K^3=0$. Hence $$\rho(W)=13.$$

\subsection{Parallelisable webs :}

These are webs such that all $u_i$'s   are affine  functions relatively to some system of local coordinates. Then, with these coordinates,  the only blocks  $M^{(h)}_k$ which are not zero in the matrices ${\cal M}_k$ are the diagonal blocks $P_k=M^{(k-1)}_k$, and   $Rank ({\cal M}_k)= \sum_{h=1}^k Rank (P_h).$
 Thus 
$$\rho_{h+1}=\rho_h+\bigl(d-Rank (P_{h+2})\bigr).$$
\indent In particular,  if the web is ordinary,   $\rho_{h+1}=\rho_h$ for $h\geq h_0-2$.  \hb 
Hence\footnote{This has already been quoted in [CL] (theorem 6-5) by other considerations.}, {\it  all ordinary parallelisable webs have maximal rank $\pi'(n,d)\ \bigl(=\rho_{h_0-2}\bigr)$}.

If they are not ordinary, and if there exists some $h_1$ ($\geq h_0-2$) such that $\rho_{h_1+1}=\rho_{h_1}$, then the sequence of the $\rho_h$'s is stationary from there because of the lemma 2-3 above, and then 
$$\rho(W)= \rho_{h_1} \ \bigl(>\pi'(n,d)\bigr).$$Such an example is given below. 

\subsection{Non ordinary  example  $n=3,\ d= 10$ :}

Let $W_{10}$  be the parallel 10-sub-web  of the  ordinary 11-web above, obtained by deleting $u_{11}$. It is not ordinary (since $P_3$ has rank 9, not  10). 
We  get :$$\pi'(3,10)=11=\rho_3<\rho_4=\rho_5= 12=\rho(W_{10})<\pi(3,10)=16.$$




   \n Jean-Paul Dufour, former professor at the University of {\cal M}ontpellier II, \hb
  1 rue du Portalet, 34820 Teyran, France\hb  email : dufourh@netcourrier.com,

  \n Daniel Lehmann, former professor at the University of {\cal M}ontpellier II,\hb   4 rue Becagrun,  30980 Saint Dionisy, France\hb  email : lehm.dan@gmail.com,

 \end{document}